\documentclass[11pt]{article}

 \usepackage{a4}                          
 \usepackage{amsmath}
 \usepackage{amssymb}
 \usepackage{graphicx}

 \newtheorem{thm}{Theorem}[section]

 \newtheorem{cor}[thm]{Corollary}

 \newtheorem{lemma}[thm]{Lemma}
 \newtheorem{rem}[thm]{Remark}

\newcommand{\Int}{\int_{-\infty}^\infty}

 \title{Linearization coefficients of Bessel polynomials}
 \author{Christian Berg and Christophe Vignat}

 \date{\today}

\begin{document}

 \maketitle

 \begin{abstract} 
We prove positivity results about linearization 
and connection coefficients for Bessel polynomials.
The proof is based on a recursion formula and explicit formulas for
the coefficients in special cases. The result implies that the
distribution of a convex combination of independent Student-t random
variables with arbitrary odd degrees of freedom has a density which
is a convex combination of certain Student-t densities with odd
degrees of freedom.
\end{abstract}

\noindent 
2000 {\em Mathematics Subject Classification}:\\
primary 33C10; secondary 60E05 

\noindent
Keywords: Bessel polynomials, Student-t distribution, linearization 
coefficients

\section{Introduction}

In this paper we consider the Bessel polynomials $q_n$ of degree $n$ 
\begin{equation}\label{eq:qn}
q_{n}\left(u\right)=\sum_{k=0}^{n}\alpha_{k}^{(n)}u^{k},
\end{equation}
where 
\begin{equation}\label{eq:alpha}
\alpha_{k}^{(n)}=\tfrac{\binom{n}{k}}{\binom{2n}{k}}\tfrac{2^k}{k!}=
\frac{n!\,(2n-k)!\,2^k}{(2n)!\,(n-k)!\,k!}.
 \end{equation}
The first examples of these polynomials
are
\[
q_{0}\left(u\right)=1, q_{1}\left(u\right)=1+u, q_{2}\left(u\right)=1+u+\frac{u^{2}}{3}.
\]
They are normalized according to 
\[
q_{n}\left(0\right)=1,
\]
and thus differ from the polynomials $\theta_{n}\left(u\right)$ in 
\cite{Grosswald}
by the constant factor $\tfrac{\left(2n\right)!}{n!2^{n}},$ i.e.
$$
\theta_n(u)=\frac{(2n)!}{n!2^{n}}q_n(u).
$$
For $\nu>0$ we recall that the probability density on $\mathbb R$
\begin{equation}\label{eq:student}
f_{\nu}(x)=\frac{A_{\nu}}{(1+x^2)^{\nu+\tfrac12}},\quad A_{\nu}=
\frac{\Gamma(\nu+\tfrac12)}{\Gamma(\tfrac12)\Gamma(\nu)}
\end{equation}
has the characteristic function
\begin{equation}\label{eq:char}
\Int e^{ixy}f_{\nu}(x)\,dx=k_{\nu}(|y|),\quad y\in\mathbb R,
\end{equation}
where
\begin{equation}\label{eq:Bes1}
k_{\nu}(u)=\frac{2^{1-\nu}}{\Gamma\left(\nu\right)}u^{\nu}K_{\nu}\left(u\right),\quad
u\geq 0,
\end{equation}
and $K_{\nu}$ is the modified Bessel function of the third kind. If 
$\nu=n+\tfrac12$ with $n=0,1,2,\ldots$ then
\begin{equation}\label{eq:Bes2}
k_{\nu}(u)=e^{-u}q_n(u),\quad u\geq 0,
\end{equation}
and $f_{\nu}$ is called a Student-t density with $2\nu=2n+1$ degrees of
freedom. For $\nu=\tfrac12$ then $f_{\nu}$ is density of a Cauchy
distribution. Note that for simplicity we have avoided the usual
scaling of the Student-t distribution.

In this paper, we provide the solutions of the three following problems:
assuming in the rest of the paper that $a$ is a constant with $0\le a\le1$,

\begin{enumerate}
\item positivity and explicit values of the connection coefficients
 $c_{k}^{\left(n\right)}(a)$
in the expansion
\begin{equation}\label{eq:pb1}
q_{n}\left(au\right)=\sum_{k=0}^{n}c_{k}^{\left(n\right)}\left(a\right)q_{k}
\left(u\right),
\end{equation}

\item positivity and explicit value of the linearization coefficients
 $\beta_{i}^{\left(n\right)}(a)$ in
the expansion\begin{equation}
q_{n}(au)q_{n}((1-a)u)=\sum_{i=0}^{n}\beta_{i}^{\left(n\right)}(a)q_{n+i}
\left(u\right),\label{eq:pb2}\end{equation}

\item positivity of the linearization coefficients
 $\beta_{k}^{\left(n,m\right)}(a)$ in
the expansion\begin{equation}
q_{n}(au)q_{m}((1-a)u)=\sum_{k=n\wedge m}^{n+m}\beta_{k}^{\left(n,m\right)}(a)
q_{k}\left(u\right).\label{eq:pb3}\end{equation}

\end{enumerate}
Note that $\beta_i^{(n)}(a)=\beta_{n+i}^{(n,n)}(a)$ and that 
(\ref{eq:pb1}) is a special case of (\ref{eq:pb3})
corresponding to $m=0$ with $c_k^{(n)}(a)=\beta_k^{(n,0)}(a)$.
 Note also that $u=0$ in (\ref{eq:pb3})
yields 
$$
\sum_{k=n\wedge m}^{n+m}\beta_k^{(n,m)}(a)=1,
$$
so (\ref{eq:pb3}) is a convex combination. As polynomial identities, 
(\ref{eq:pb1})-(\ref{eq:pb3}) of course hold for all complex $a,u$,
but as we will see later, the positivity of the coefficients holds
only for $0\leq a\leq 1.$ 

Because of (\ref{eq:char}) and (\ref{eq:Bes2})
formula (\ref{eq:pb3}) is equivalent with the following identity
between Student-t densities 
\begin{equation}\label{eq:conv}
\frac{1}{a}f_{n+\tfrac12}\left(\frac{x}{a}\right)*\frac{1}{1-a}
f_{m+\tfrac12}\left(\frac{x}{1-a}\right)=\sum_{k=n\wedge m}^{n+m}\beta_k^{(n,m)}(a)f_{k+\tfrac12}(x)
\end{equation}
for $0<a<1$ and $*$ is the ordinary convolution of densities.

Although (\ref{eq:pb3}) is more general than
(\ref{eq:pb1}),(\ref{eq:pb2}),
we stress that we give explicit formulas below for $c_k^{(n)}(a)$ and
$\beta_i^{(n)}(a)$ from which the positivity is clear. The positivity
of $\beta_k^{(n,m)}(a)$ for the general case can be deduced from the
special cases via a recursion formula, see Lemma \ref{thm:betarec} below.

These problems have an important application in statistics: the Behrens-Fisher
problem consists in testing the equality of the means of two normal
populations. Fisher \cite{Fisher}%
\footnote{the collected papers of R.A. Fisher are available at the following
address

http://www.library.adelaide.edu.au/digitised/fisher/%
} has shown that this test can be performed using the $d-$statistics
defined as \[
d_{f_{1},f_{2},\theta}=t_{1}\sin\theta-t_{2}\cos\theta,\]
where $t_{1}$ and $t_{2}$ are two independent Student-t random variables
with respective degrees of freedom $f_{1}$ and $f_{2}$ and 
$\theta\in\left[0,\tfrac{\pi}{2}\right]$. Many different
results have been obtained on the behaviour of the $d-$statistics.
Tables of the distribution of $d_{f_{1},f_{2},\theta}$ have been
provided in 1938 by Sukhatme \cite{sukhatme} at Fisher's suggestion.
In 1956, Fisher and Healy explicited the distribution of $d_{f_{1},f_{2},\theta}$
as a mixture of Student-t distributions (Student-t distribution with
a random, discrete number of degrees of freedom) for small, odd values
of $f_{1}$ and $f_{2}$. This work was extended by Walker and Saw
\cite{walker} who provided, still in the case of odd numbers of degrees
of freedom, an explicit way of computing the coefficients of the Student-t
mixture as solutions of a linear system; however, they did not prove
the positivity of these coefficients, claiming only 

\begin{quotation}
{}``Extensive numerical investigation indicates also that $\eta_{i}\ge0$
for all i; however, an analytic proof has not been found.''
\end{quotation}
This conjecture is proved in Theorem 2 and 3 below. Section 2 of this
paper gives the explicit solutions to problems 1, 2 and 3, whereas
section 3 is dedicated to their proofs. The last section gives an extension
of Theorem 2 in terms of inverse Gamma distributions.

Using the fact that the Student-t distribution is a scale mixture of normal
distributions by an inverse gamma distribution our positivity result
is equivalent to an analogous positivity result for inverse gamma
distributions. This result has been observed for small values of 
the degrees of freedom in
\cite{witkowsky}. In \cite{Giron} the coefficients are claimed to be
non-negative but the paper does not contain any arguments to prove it.

\section{Results}

\subsection{Solution of problem 1 and a stochastic interpretation}

\begin{thm}
\label{thm:1}The coefficients $c_{k}^{\left(n\right)}\left(a\right)$
in (\ref{eq:pb1}) write\[
c_{k}^{\left(n\right)}\left(a\right)=a^{k}\tfrac{\binom{n}{k}}{\binom{2n}{2k}}\sum_{r=1}^{\left(n-k\right)\wedge\left(k+1\right)}\binom{n+1}{k+1-r}\binom{n-k-1}{r-1}\left(1-a\right)^{r}\]
for $0\leq k\leq n-1$ while $c_n^{(n)}(a)=a^n$ and, as $0\le a\le1,$
they are positive.
\end{thm}
A stochastic interpretation of Theorem \ref{thm:1} writes as follows:
replacing u by $\vert u\vert$ and multiplying equation (\ref{eq:pb1})
by $\exp\left(-\vert u\vert\right)$, we obtain
 \begin{equation}
e^{-\left(1-a\right)\vert u\vert}e^{-a|u|}q_{n}\left(a\vert
  u\vert\right)=
\sum_{k=0}^{n}c_{k}^{\left(n\right)}\left(a\right)q_{k}(|u|)e^{-|u|}.
\label{eq:interp1}\end{equation}
Equality (\ref{eq:interp1}) can be interpreted as follows: the convex
combination of an independent Cauchy variable $C$ and a Student-t
variable $X_{n}$ with $2n+1$ degrees of freedom follows a
 Student-t distribution with
random number $2K\left(\omega\right)+1$ of degrees of freedom:
\[
\left(1-a\right)C+aX_{n}\overset{d}{=}X_{K\left(\omega\right)},\]
where $K\left(w\right)\in[0,n]$ is a discrete random variable such
that \[
\Pr\{ K\left(w\right)=k\}=c_{k}^{\left(n\right)}\left(a\right),\quad 0\le k\le n.\]

\subsection{Solution of problem 2 and a probabilistic interpretation}

\begin{thm}
\label{thm:2}The coefficients $\beta_{i}^{\left(n\right)}\left(a\right)$
in (\ref{eq:pb2}) write
\begin{eqnarray*}
\beta_{i}^{(n)}(a)&=&(4a(1-a))^{i}\left(\frac{n!}{(2n)!}\right)^{2}2^{-2n}
\frac{(2n-2i)!(2n+2i)!}{(n-i)!(n+i)!}\\
&\times&\sum_{j=0}^{n-i}
\binom{2n+1}{2j}\binom{n-j}{i}(2a-1)^{2j}
\end{eqnarray*}
and, as $0\le a\le1,$ they are positive.
\end{thm}

A probabilistic interpretation of this result writes as follows.

\begin{cor}
With $a=\sin\theta$, $f_{1}=f_{2}=2n+1,$ statistic $d_{f_{1},f_{2},\theta}$
follows a Student-t distribution with a random number of degrees of
freedom $F\left(\omega\right)$ distributed according to\[
\Pr\left\{ F\left(\omega\right)=2n+2i+1\right\}
=\beta_{i}^{\left(n\right)}\left(a\right),\quad 0\leq i\leq n.\]
 \end{cor}

\subsection{Problem 3}

\begin{thm}
\label{thm:3} The coefficients $\beta_{k}^{\left(n,m\right)}(a)$ in (\ref{eq:pb3})
are positive for $0\le a\le1.$
\end{thm}
We were unable to derive the explicit values of the coefficients 
$\beta_{k}^{\left(n,m\right)}(a)$:
however, their positivity allows us to claim the following Corollary.

\begin{cor}
With $a=\sin\theta$, $f_{1}=2n+1, f_{2}=2m+1,$ statistic $d_{f_{1},f_{2},\theta}$
follows a Student-t distribution with a random number of degrees of
freedom $F\left(\omega\right)$ distributed according to\[
\Pr\left\{ F\left(\omega\right)=2k+1\right\}
=\beta_{k}^{\left(n,m\right)}\left(a\right),\quad n\wedge m\leq k\leq n+m.\]
\end{cor}

\begin{thm}\label{thm:4}
For $k\geq 2$ let $n_1,\ldots,n_k$ be nonnegative integers and let
$a_1,\ldots,a_k$ be positive real numbers with sum 1. Then 
\begin{equation}\label{eq:final}
q_{n_1}(a_1u)q_{n_2}(a_2u)\cdots
q_{n_k}(a_ku)=\sum_{j=l}^L\beta_jq_j(u),\quad u\in\mathbb R
\end{equation}
with nonnegative coefficients $\beta_j$ with sum 1 and 
$l=\min(n_1,\ldots,n_k),L=n_1+\cdots +n_k$.
\end{thm}

\section{Proofs}
\subsection{Generalities about Bessel polynomials}

As a preparation to the proofs we give some recursion formulas for
$q_n$. They follow from corresponding formulas for $\theta_n$ from 
\cite{Grosswald}, but they can also be proved directly from the
definitions (\ref{eq:qn}) and (\ref{eq:alpha}).
The formulas are
\begin{equation}\label{eq:rec1}
q_{n+1}(u)=q_n(u)+\frac{u^2}{4n^2-1}q_{n-1}(u),\quad n\geq 1,
\end{equation}

\begin{equation}\label{eq:dif}
q_{n}'(u)=q_n(u)-\frac{u}{2n-1}q_{n-1}(u),\quad n\geq 1.
\end{equation}
We can write 
\begin{equation}\label{eq:inverse}
u^n=\sum_{i=0}^n \delta_i^{(n)}q_i(u),\quad n=0,1,\ldots
\end{equation}
and $\delta_i^{(n)}$ is given by a formula due to Carlitz \cite{C},
see \cite[p. 73]{Grosswald} or \cite{walker}:
\begin{equation}\label{eq:carlitz}
 \delta_{i}^{\left(n\right)}=\left\{ \begin{array}{cc}
\frac{\left(n+1\right)!}{2^{n}}\frac{\left(-1\right)^{n-i}\left(2i\right)!}{\left(n-i\right)!i!\left(2i+1-n\right)!} & \text{for }\frac{n-1}{2}\le i\le n\\
0 & \text{for } 0\le i<\frac{n-1}{2}\end{array}\right..
\end{equation}

Later we need the following extension of (\ref{eq:rec1}) which we
formulate using the Pochhammer symbol $(z)_n:=z(z+1)\cdots (z+n-1)$
for $z\in\mathbb C, n=0,1,\ldots.$

\begin{lemma}\label{thm:extrec}
For $0\leq k\leq n$ we have \[
u^{2k}q_{n-k}(u)=\sum_{i=0}^{k}\gamma^{(n,k)}_{i}q_{n+i}(u)\]
 where \begin{equation}
\gamma^{(n,k)}_{i}=2^{2k}\binom{k}{i}(n-k+\tfrac{1}{2})_{k+i}(-n-\tfrac{1}{2})_{k-i}.\label{eq:coef1}\end{equation}
\end{lemma}

{\it Proof:}
The Lemma is trivial for $k=0$ and reduces to the recursion (\ref{eq:rec1})
for $k=1$ written as
\begin{equation}
u^{2}q_{n-1}(u)=2^2(n-\tfrac{1}{2})_{2}\left(q_{n+1}(u)-q_{n}(u)\right).
\label{eq:recursion}\end{equation}
 We will prove the formula (\ref{eq:coef1}) by induction in $n$,
so assume it holds for some $n$ and all $0\leq k\leq n$. Multiplying
the formula of the lemma by $u^{2}$ we get \[
u^{2k+2}q_{n-k}(u)=\sum_{i=0}^{k}\gamma^{(n,k)}_{i}u^{2}q_{n+i}(u),\]
 hence by (\ref{eq:recursion}) 
\begin{eqnarray*}
 &  & u^{2(k+1)}q_{n+1-(k+1)}(u)=\sum_{i=0}^{k}\gamma^{(n,k)}_{i}2^{2}
(n+i+\tfrac{1}{2})_2\left[q_{n+i+2}(u)-q_{n+i+1}(u)\right]\\
 & = & \gamma^{(n,k)}_{k}2^{2}(n+k+\tfrac{1}{2})_2\,q_{n+k+2}(u)\\
 & + & \sum_{i=1}^{k}2^{2}(n+i+\tfrac{1}{2})\left[\gamma^{(n,k)}_{i-1}
(n+i-\tfrac{1}{2})-\gamma^{(n,k)}_{i}(n+i+\tfrac{3}{2})\right]q_{n+1+i}(u)\\
 & - &
 \gamma^{(n,k)}_{0}2^{2}(n+\tfrac{1}{2})_2\,q_{n+1}(u).
\end{eqnarray*}
 Using the induction hypothesis we easily get \[
\gamma^{(n,k)}_{k}2^{2}(n+k+\tfrac{1}{2})_2=2^{2k+2}(n-k+
\tfrac{1}{2})_{2k+2}=\gamma^{(n+1,k+1)}_{k+1},\]
 and \[
-\gamma^{(n,k)}_{0}2^{2}(n+\tfrac{3}{2})(n+\tfrac{1}{2})=2^{2k+2}(n-k+\tfrac{1}{2})_{k+1}(-n-\tfrac{3}{2})_{k+1}=\gamma^{(n+1,k+1)}_{0}.\]
 Concerning the coefficient $C$ to $q_{n+1+i}(u)$ above we have
\[
C=2^{2k+2}(n+i+\tfrac{1}{2})\left[\binom{k}{i-1}(n-k+\tfrac{1}{2})_{k+i-1}
(-n-\tfrac{1}{2})_{k-i+1}(n+i-\tfrac{1}{2})\right.\]
 \[
-\left.\binom{k}{i}(n-k+\tfrac{1}{2})_{k+i}(-n-\tfrac{1}{2})_{k-i}(n+i+\tfrac{3}{2})\right]\]
\begin{eqnarray*}
 & = & 2^{2k+2}(n-k+\tfrac{1}{2})_{k+1+i}(-n-\tfrac{1}{2})_{k-i}\\
& \times &\left[\binom{k}{i-1}(-n-\tfrac{1}{2}+k-i)-\binom{k}{i}(n+i+\tfrac{3}{2})\right]\\
 & = & 2^{2k+2}(n-k+\tfrac{1}{2})_{k+1+i}(-n-\tfrac{1}{2})_{k-i}\left[\binom{k+1}{i}(-n-\tfrac{3}{2})\right]\\
 & = & 2^{2k+2}\binom{k+1}{i}(n-k+\tfrac{1}{2})_{k+1+i}(-n-\tfrac{3}{2})_{k+1-i}=\gamma^{(n+1,k+1)}_{i}.\end{eqnarray*}
$\square$

We stress that Lemma \ref{thm:extrec} is the special case
$\nu=n+\tfrac12$ of the
following recursion for modified Bessel functions of the third kind.

\begin{lemma}
For all $\nu>0$ and all nonnegative integers $j<\nu$ we have for $u>0$
\[
u^{\nu+j}K_{\nu-j}\left(u\right)=\sum_{i=0}^{j}(-2)^{j-i}\binom{j}{i}
\frac{\Gamma(\nu+1)}{\Gamma(\nu+1-(j-i))}
u^{\nu+i}K_{\nu+i}\left(u\right)\]
and
\[
u^{2j}k_{\nu-j}\left(u\right)=\sum_{i=0}^{j}\left(-1\right)^{j-i}2^{2j}\binom{j}{i}\frac{\Gamma\left(\nu+1\right)\Gamma\left(\nu+i\right)}{\Gamma\left(\nu+1-\left(j-i\right)\right)\Gamma\left(\nu-j\right)}k_{\nu+i}\left(u\right).\]
\end{lemma}

{\it Proof:} The second formula follows from the first using formula
 (\ref{eq:Bes1}), and the first can be proved by induction using the
 following recursion formula for modified Bessel functions of the
 third kind, cf. \cite[p. 79]{Wa}
$$ 
K_{\nu-1}(u)=K_{\nu+1}(u)-\frac{2\nu}{u}K_\nu(u).
$$
We skip the details. $\square$

\subsection{Proof of Theorem \ref{thm:1}}

From (\ref{eq:qn}) and (\ref{eq:inverse}) we get
\[
q_{n}\left(au\right)=\sum_{j=0}^{n}\alpha_{j}^{\left(n\right)}a^{j}
\sum_{i=0}^{j}\delta_{i}^{\left(j\right)}q_{i}\left(u\right)=
\sum_{k=0}^{n}c_{k}^{\left(n\right)}(a)q_{k}\left(u\right)\]
with 
\begin{eqnarray*}
c_{k}^{\left(n\right)}\left(a\right)&=&\sum_{j=k}^{n}a^{j}
\alpha_{j}^{\left(n\right)}\delta_{k}^{\left(j\right)}\\
&=&a^{k}\frac{n!}{\left(2n\right)!}\frac{\left(2k\right)!}{k!}
\sum_{j=k,j\le2k+1}^{n}\left(-a\right)^{j-k}
\frac{\left(2n-j\right)!\left(j+1\right)}{\left(n-j\right)!\left(j-k\right)!
\left(2k+1-j\right)!}
\end{eqnarray*}

In particular $c_{n}^{\left(n\right)}\left(a\right)=a^{n}$ and for
$0\le k\le n-1$
\begin{equation}\label{eq:ck}
c_{k}^{\left(n\right)}\left(a\right)=a^{k}\frac{n!}{\left(2n\right)!}
\frac{\left(2k\right)!}{k!}p(a),
\end{equation}
where
$$
p(a)=\sum_{i=0}^{(n-k)\wedge(k+1)}
(-a)^{i}\frac{(2n-k-i)!(k+i+1)}{(n-k-i)!i!(k+1-i)!}.
$$
We clearly have
$$
p(a)=\sum_{r=0}^{(n-k)\wedge(k+1)}(-1)^r\frac{p^{(r)}(1)}{r!}(1-a)^r
$$
with
$$
p^{(r)}(1)=\sum_{i=r}^{(n-k)\wedge(k+1)}(-1)^{i}\frac{(2n-k-i)!(k+i+1)}
{(n-k-i)!(i-r)!(k+1-i)!}
$$
and we only consider $0\leq r\leq (n-k)\wedge(k+1)$. To sum this
we shift the summation by $r$. For simplicity we define
$T:=(n-k-r)\wedge(k+1-r)$ and get
$$
(-1)^{r}p^{(r)}(1)=\sum_{i=0}^T (-1)^{i}\frac{(2n-k-r-i)!
\,(k+r+i+1)}{(n-k-r-i)!\,i!\,(k+1-r-i)!}.
$$
We write $k+r+1+i=(2k+2)-(k+1-r-i)$ and split the above sum
accordingly
\begin{eqnarray*}
(-1)^{r}p^{(r)}(1)&=&(2k+2)\sum_{i=0}^T (-1)^{i}
\frac{(2n-k-r-i)!}{(n-k-r-i)!\,i!\,(k+1-r-i)!}\\
&-&\sum_{i=0}^T (-1)^{i}\frac{(2n-k-r-i)!}
{(n-k-r-i)!\,i!\,(k-r-i)!}.
\end{eqnarray*}
Note that for nonnegative integers $a,b,c$ with $b,c\leq a$ we have
$$
\sum_{i=0}^{b\wedge
  c}(-1)^{i}\frac{(a-i)!}{(b-i)!(c-i)!i!}=\frac{a!}{b!c!}
\sum_{i=0}^{b\wedge c}
\frac{(-b)_i(-c)_i}{(-a)_i i!}=\frac{a!}{b!c!}\,{}_2F_1(-b,-c;-a;1),
$$
where we use that the
 sum is  an ${}_2F_1$
evaluated at 1. Its value is given by the Chu-Vandermonde formula, see
\cite{A:A:R},
 hence 
$$
\sum_{i=0}^{b\wedge c}(-1)^{i}\frac{(a-i)!}{(b-i)!(c-i)!i!}
=\frac{a!(c-a)_b}{(-a)_b
  b!c!}.
$$
The two sums above are of this form and we get
$$
(-1)^{r}p^{(r)}(1)=\frac{(2n-k-r)!}{(n-k-r)!(k+1-r)!}Q,
$$
where
$$
Q=(2k+2)\frac{(2k-2n+1)_{n-k-r}}{(k+r-2n)_{n-k-r}}-
(k+1-r)\frac{(2k-2n)_{n-k-r}}{(k+r-2n)_{n-k-r}}
$$
$$
=\frac{(2k-2n+1)_{n-k-r-1}}{(k+r-2n)_{n-k-r}}[(2k+2)(k-r-n)-(k+1-r)(2k-2n)]
$$
$$
=2r(n+1)\frac{(n+r+1-k)_{n-k-r-1}}{(n+1)_{n-k-r}},
$$
where we used $(a)_n=(-1)^{n}(1-a-n)_n$ twice.
This gives
$$
(-1)^{r}p^{(r)}(1)=2r\binom{n+1}{k+1-r}\frac{(2n-2k-1)!}{(n-k-r)!}
 $$
and finally
$$
p(a)=\sum_{r=0}^{(n-k)\wedge(k+1)} (1-a)^r\frac{2r}{r!}
\binom{n+1}{k+1-r}\frac{(2n-2k-1)!}{(n-k-r)!}.
$$
Note that the term corresponding to $r=0$ is zero. If we insert this expression
for $p(a)$ in (\ref{eq:ck}), we get the formula of Theorem \ref{thm:1}.
$\square$

\begin{rem} {\rm The evaluation above of $(-1)^rp^{(r)}(1)$ can be
    done using generating functions like in \cite{walker}. The authors
    want to thank Mogens Esrom Larsen for the idea to use the
    Chu-Vandermonde identity twice.}
\end{rem}
 
\subsection{Proof of Theorem \ref{thm:2}}

The starting point is the following formula of Macdonald, see \cite{Wa}\begin{equation}
K_{\nu}(z)K_{\nu}(X)=\frac{1}{2}\int_{0}^{\infty}\exp[-\frac{s}{2}-\frac{z^{2}+X^{2}}{2s}]K_{\nu}(\frac{zX}{s})\frac{ds}{s},\label{eq:Macdonald}\end{equation}
 which we will use for $\nu=n+\frac{1}{2}, z=au, X=(1-a)u$. Multiplying
(\ref{eq:Macdonald}) by \[
\left(\frac{2^{1-\nu}}{\Gamma(\nu)}\right)^{2}(a(1-a)u^{2})^{\nu}\]
 and using (\ref{eq:Bes1})
 we find 
$$
k_{\nu}(au)k_{\nu}((1-a)u)=
$$
$$
\frac{1}{2^{\nu}\Gamma(\nu)}\int_{0}^{\infty}\exp\left[-\frac{s}{2}-u^{2}\frac{a^{2}+(1-a)^{2}}{2s}\right]s^{\nu-1}k_{\nu}\left(\frac{a(1-a)u^{2}}{s}\right)\;ds.
$$
 We now insert that with $\nu=n+\frac{1}{2}$ we have
$k_{\nu}(|u|)=e^{-|u|}q_{n}(|u|)$
 and hence after some simplification
$$
e^{-|u|}q_{n}(a|u|)q_{n}((1-a)|u|)=
$$
$$
\frac{1}{2^{n+\frac{1}{2}}\Gamma(n+\frac{1}{2})}\int_{0}^{\infty}\exp\left[-\frac{s}{2}-\frac{u^{2}}{2s}\right]s^{n-\frac{1}{2}}q_{n}\left(\frac{a(1-a)u^{2}}{s}\right)\;ds.
$$
We next insert the expression (\ref{eq:qn}) for $q_{n}$ under the
integral sign. This
gives \[
e^{-|u|}q_{n}(a|u|)q_{n}((1-a)|u|)=\]
 \[
\sum_{k=0}^{n}\alpha_{k}^{(n)}(a(1-a))^{k}u^{2k}\frac{1}{2^{n+\frac{1}{2}}\Gamma(n+\frac{1}{2})}\int_{0}^{\infty}\exp\left[-\frac{s}{2}-\frac{u^{2}}{2s}\right]s^{n-k-\frac{1}{2}}\; ds.\]
 Using the following formula from \cite[3.471(9)]{G:R} 
\begin{equation}\label{eq:GR}
\int_{0}^{\infty}x^{\nu-1}\exp(-\frac{\beta}{x}-\gamma x)\;
dx=2\left(\frac{\beta}{\gamma}\right)^{\nu/2}K_{\nu}\left(2\sqrt{\beta\gamma}\right)
\end{equation}
and again (\ref{eq:Bes1}) the above is equal to
 \[
=\sum_{k=0}^{n}\alpha_{k}^{(n)}(a(1-a))^{k}\frac{2^{n-k+\frac{1}{2}}\Gamma(n-k+\frac{1}{2})}{2^{n+\frac{1}{2}}\Gamma(n+\frac{1}{2})}e^{-|u|}u^{2k}q_{n-k}(|u|).\]
Finally, using Pochhammer symbols and skipping absolute values
since we are now dealing with a polynomial identity, we get
 \begin{equation}
q_{n}(au)q_{n}((1-a)u)=\sum_{k=0}^{n}\alpha_{k}^{(n)}(a(1-a))^{k}\frac{(\frac{1}{2})_{n-k}}{2^{k}(\frac{1}{2})_{n}}u^{2k}q_{n-k}(u).\label{eq:Macdonald3}\end{equation}

Using the expression for $u^{2k}q_{n-k}(u)$ from Lemma  \ref{thm:extrec}
and the expression for $\alpha_{k}^{(n)}$ in (\ref{eq:Macdonald3})
we then get \begin{eqnarray*}
q_{n}(au)q_{n}((1-a)u) & = &
\sum_{k=0}^{n}\frac{\binom{n}{k}(\frac{1}{2})_{n-k}}{\binom{2n}{k}(\frac{1}{2})_{n}k!}(a(1-a))^{k}\sum_{i=0}^{k}
\gamma^{(n,k)}_{i}q_{n+i}(u)\\
 & = &
 \sum_{i=0}^{n}q_{n+i}(u)\sum_{k=i}^{n}(a(1-a))^{k}\frac{\binom{n}{k}(\frac{1}{2})_{n-k}}{\binom{2n}{k}(\frac{1}{2})_{n}k!}\gamma^{(n,k)}_{i}
\end{eqnarray*}
 hence \[
q_{n}(au)q_{n}((1-a)u)=\sum_{i=0}^{n}\beta_{i}^{(n)}(a)q_{n+i}(u)\]
 with \begin{eqnarray*}
\beta_{i}^{(n)}(a) & = & \sum_{k=i}^{n}(a(1-a))^{k}\frac{\binom{n}{k}(\frac{1}{2})_{n-k}}{\binom{2n}{k}(\frac{1}{2})_{n}k!}2^{2k}\binom{k}{i}(n-k+\frac{1}{2})_{k+i}(-n-\frac{1}{2})_{k-i}\\
 & = & (a(1-a))^{i}\sum_{l=0}^{n-i}(a(1-a))^{l}\frac{\binom{n}{i+l}(\frac{1}{2})_{n-i-l}}{\binom{2n}{i+l}(\frac{1}{2})_{n}(i+l)!}2^{2i+2l}\binom{i+l}{i}\\
 & \times & (n-i-l+\frac{1}{2})_{l+2i}(-n-\frac{1}{2})_{l}\end{eqnarray*}
 Collecting \[
(\frac{1}{2})_{n-i-l}(n-i-l+\frac{1}{2})_{l+2i}=(\frac{1}{2})_{n+i}\]
 we get 
$$
\beta_{i}^{(n)}(a)=(a(1-a))^{i}\sum_{l=0}^{n-i}(a(1-a))^{l}
\frac{\binom{n}{i+l}(\frac{1}{2})_{n+i}}{\binom{2n}{i+l}(\frac{1}{2})_{n}}
\frac{2^{2i+2l}}{i!\,l!}(-n-\frac{1}{2})_{l}
$$
$$
  =(a(1-a))^{i}\frac{n!(\frac{1}{2})_{n+i}2^{2i}}{(2n)!\,(\frac{1}{2})_{n}\,
i!}\sum_{l=0}^{n-i}(4a(1-a))^{l}\frac{(2n-i-l)!\,(-n-\frac{1}{2})_{l}}
{(n-i-l)!\,l!}
$$
$$
  = (a(1-a))^{i}\left(\frac{n!}{(2n)!}\right)^{2}
\frac{(2n+2i)!}{(n+i)!i!}\sum_{l=0}^{n-i}\left(1-(2a-1)^{2}\right)^{l}
\frac{(2n-i-l)!\,(-n-\frac{1}{2})_{l}}{(n-i-l)!\,l!}.
$$
Expanding $(1-(2a-1)^{2})^{l}$ using the binomial formula and interchanging
the sums we get \[
\sum_{l=0}^{n-i}(1-(2a-1)^{2})^{l}\frac{(2n-i-l)!(-n-\frac{1}{2})_{l}}{(n-i-l)!l!}\]
 \[
=\sum_{j=0}^{n-i}(-1)^{j}\frac{(2a-1)^{2j}}{j!}\sum_{l=j}^{n-i}\frac{(2n-i-l)!(-n-\frac{1}{2})_{l}}{(n-i-l)!(l-j)!}.\]

We claim that \begin{eqnarray}
S:&=&\sum_{l=j}^{n-i}\frac{(2n-i-l)!(-n-\frac{1}{2})_{l}}{(n-i-l)!(l-j)!}
\nonumber\\
&=&\frac{(2n-2i)!}{(n-i)!}2^{2i-2n}(-1)^{j}j!i!\binom{2n+1}{2j}\binom{n-j}{i}.\label{eq:CV}\end{eqnarray}
 and we have obtained the final formula 
\begin{eqnarray*}
\beta_{i}^{(n)}(a)=(4a(1-a))^{i}\left(\frac{n!}{(2n)!}\right)^{2}2^{-2n}
\frac{(2n-2i)!(2n+2i)!}{(n-i)!(n+i)!}\\
\times\sum_{j=0}^{n-i}\binom{2n+1}{2j}\binom{n-j}{i}(2a-1)^{2j}.
\end{eqnarray*}
 We will see that (\ref{eq:CV}) is a Chu-Vandermonde formula. In
fact, shifting the summation index putting $l=j+m$ we get \[
S=(-n-\frac{1}{2})_{j}\sum_{m=0}^{n-i-j}\frac{(2n-i-j-m)!(-n+j-\frac{1}{2})_{m}}{m!(n-i-j-m)!},\]
 and calling the general term in this sum $c_{m}$ we get \[
\frac{c_{m+1}}{c_{m}}=\frac{(m+i+j-n)(m+j-n-\frac{1}{2})}{(m+1)(m+i+j-2n)},\]
 which shows that the sum is a $_{2}F_{1}$. We have \[
S=(-n-\tfrac12)_j\,c_0\times\,{}_{2}F_{1}(-(n-i-j),j-n-\frac{1}{2};i+j-2n;1),\]
 and using the Chu-Vandermonde identity, cf. \cite{A:A:R}:
 \[
_{2}F_{1}(-n,a;c;1)=\frac{(c-a)_{n}}{(c)_{n}},\]
 we get \begin{eqnarray*}
S & = & (-n-\frac{1}{2})_{j}\frac{(2n-i-j)!}{(n-i-j)!}\frac{(-n+i+\frac{1}{2})_{n-i-j}}{(i+j-2n)_{n-i-j}}\\
 & = & (-n-\frac{1}{2})_{j}\frac{(2n-i-j)!}{(n-i-j)!}\frac{(j+\frac{1}{2})_{n-i-j}}{(n+1)_{n-i-j}},\end{eqnarray*}
 where we used $(a)_n=(-1)^n(1-a-n)_n$ twice.
We can now simplify to get \[
S=(-n-\frac{1}{2})_{j}\frac{n!}{(n-i-j)!}\frac{(\frac{1}{2})_{n-i}}{(\frac{1}{2})_{j}},\]
 and applying the formula \[
(\frac{1}{2})_{p}=\frac{(2p)!}{p!2^{2p}}\]
 twice we get \[
S=\frac{(2n-2i)!}{(n-i)!}2^{2i-2n}(-n-\frac{1}{2})_{j}\frac{n!}{(n-i-j)!}\frac{j!2^{2j}}{(2j)!}.\]
 Now we can write \[
(-n-\frac{1}{2})_{j}2^{2j}=(-1)^{j}\frac{(2n+1)!(n-j)!}{n!(2n-2j+1)!},\]
 and hence \[
S=\frac{(2n-2i)!}{(n-i)!}2^{2i-2n}(-1)^{j}j!i!\binom{2n+1}{2j}\binom{n-j}{i}.\]

$\square$

\subsection{Proof of Theorem \ref{thm:3}}

For $n,m\ge0$ and $a\in\mathbb{R}$, we can write 
\begin{equation}
q_{n}\left(au\right)q_{m}\left(\left(1-a\right)u\right)=\sum_{k=0}^{m+n}\beta_{k}^{\left(n,m\right)}\left(a\right)q_{k}\left(u\right)\label{eq:qnqm}\end{equation}
for some uniquely determined coefficients since the left-hand side
is a polynomial in $u$ of degree $\le n+m$. Clearly
$\beta_k^{(n,m)}(a)$ is a polynomial in $a$ satisfying
\begin{equation}
\beta_{k}^{\left(n,m\right)}\left(a\right)=\beta_k^{\left(m,n\right)}
\left(1-a\right).\label{eq:symmetry}\end{equation}
We shall prove that $\beta_k^{(n,m)}(a)\geq 0$ for $0\leq a\leq 1$ and
that $\beta_{k}^{\left(n,m\right)}\left(a\right)=0$
 if $k<n\wedge m$, which will be a consequence of the following
 recursion formula.

\begin{lemma}\label{thm:betarec}
For $n,m\ge1$, we have\begin{equation}
\frac{1}{2k+1}\beta_{k+1}^{\left(n,m\right)}\left(a\right)=\frac{a^{2}}{2n-1}\beta_{k}^{\left(n-1,m\right)}\left(a\right)+\frac{\left(1-a\right)^{2}}{2m-1}\beta_{k}^{\left(n,m-1\right)}\left(a\right),\label{eq:lemma}\end{equation}
where $k=0,1,\dots,m+n-1.$ 

Furthermore $\beta_0^{(n,m)}(a)=0.$
\end{lemma}

{\it Proof:}
Differentiating (\ref{eq:qnqm}) with respect to $u$ gives\[
aq'_{n}\left(au\right)q_{m}\left(\left(1-a\right)u\right)+\left(1-a\right)q_{n}\left(au\right)q'_{m}\left(\left(1-a\right)u\right)=\sum_{k=1}^{m+n}\beta_{k}^{\left(n,m\right)}\left(a\right)q_{k}'\left(u\right)\]
and using the formula (\ref{eq:dif})
we find\begin{eqnarray*}
 & a\left(q_{n}\left(au\right)-\frac{au}{2n-1}q_{n-1}\left(au\right)\right)q_{m}\left(\left(1-a\right)u\right)\\
+ & \left(1-a\right)q_{n}\left(au\right)\left(q_{m}\left(\left(1-a\right)u\right)-\frac{\left(1-a\right)u}{2m-1}q_{m-1}\left(\left(1-a\right)u\right)\right)\\
= &
\sum_{k=1}^{m+n}\beta_{k}^{\left(n,m\right)}\left(a\right)\left(q_{k}
\left(u\right)-\frac{u}{2k-1}q_{k-1}\left(u\right)\right)\end{eqnarray*}
and using (\ref{eq:qnqm}) once more we get
\begin{eqnarray*}
 & -\frac{a^{2}u}{2n-1}q_{n-1}\left(au\right)q_{m}\left(\left(1-a\right)u\right)-\frac{\left(1-a\right)^{2}u}{2m-1}q_{n}\left(au\right)q_{m-1}\left(\left(1-a\right)u\right)\\
= & -\beta_0^{(n,m)}(a)
-u\sum_{k=0}^{n+m-1}\beta_{k+1}^{\left(n,m\right)}
\left(a\right)(2k+1)^{-1}q_{k}(u).\end{eqnarray*}
For $u=0$ this gives $\beta_0^{(n,m)}(a)=0$ and dividing
by $-u$ and equating the coefficients of $q_{k}\left(u\right)$,
we get the desired formula.
$\square$

Now the proof of Theorem \ref{thm:3} is easy by induction in $k$ and
by the symmetry formula (\ref{eq:symmetry}),
we can assume $n\ge m.$ Let $0\leq a\leq 1$.
We  prove that $\beta_k^{(n,m)}(a)\geq 0$ for $k\leq n+m$
and that it is zero for $k<m$ (under the assumption $n\geq m$).
This is true for $k=0$ by Lemma \ref{thm:betarec} when $m\geq 1$, and
for $m=0$ it follows by Theorem \ref{thm:1}.
Assume now that the results hold for $k=k_0$ and assume $k_0+1\leq
n+m$. The nonnegativity for $k=k_0+1$ now follows by Lemma
\ref{thm:betarec}, and likewise if $k_0+1<m\leq n$ the coefficient is
0 since $k_0<(n-1)\wedge(m-1)$.
$\square$

\subsection{Proof of Theorem \ref{thm:4}}

By Theorem \ref{thm:3} the result holds for
$k=2$. Assuming it holds for $k-1\geq 2$ we have
\begin{equation}\label{eq:induction}
q_{n_1}(a_1u)\cdots
q_{n_{k-1}}(a_{k-1}u)=\sum_{j=l'}^{L'}\gamma_jq_j((1-a_k)u),
\quad u\in\mathbb R
\end{equation}
with $l'=\min(n_1,\ldots,n_{k-1}), L'=n_1+\cdots+n_{k-1}$ and
$\gamma_j\geq 0$ because we can write
$$
a_ju=\frac{a_j}{1-a_k}(1-a_k)u,\quad j=1,\ldots,k-1.
$$
If we multiply (\ref{eq:recursion}) with $q_{n_k}(a_ku)$ we get
$$
\sum_{j=l'}^{L'}\gamma_jq_{n_k}(a_ku)q_j((1-a_k)u)=
\sum_{j=l'}^{L'}\gamma_j\sum_{i=n_k\wedge
  j}^{n_k+j}\beta_i^{(n_k,j)}(a_k)q_i(u),
$$
and the assertion follows. $\square$

\section{Inverse Gamma distribution}

Grosswald proved \cite{Grosswald2} that the Student-t distribution
is infinitely divisible. This is a consequence of the infinite divisibility
of the inverse Gamma distribution because of subordination. It was
proved later that the inverse Gamma distribution is a generalized
Gamma convolution in the sense of Thorin, which is stronger than
 self-decomposability and in particular stronger than infinite divisibility. 

The following density on the half-line is an inverse Gamma density
with scale parameter $\tfrac14$ and shape parameter $\nu>0$: 

\begin{equation}\label{eq:inversegaussian}
C_{\nu}\exp(-\frac{1}{4t})t^{-\nu-1},\quad C_{\nu}=\frac{1}{2^{2\nu}
\Gamma(\nu)}.
\end{equation}

Let the corresponding probability measure be denoted
$\tilde{\gamma}_{\nu}$ and let
\[
g_{t}(x)=\frac{1}{\sqrt{4\pi t}}\exp(-\frac{x^{2}}{4t}),\quad
t>0,x\in\mathbb R
\]
 denote the Gaussian semigroup of normal densities (in the normalization
of \cite{B:F}). Then 
\begin{equation}\label{eq:subordination}
f_{\nu}(x)=\int_{0}^{\infty}g_{t}(x)\, d\tilde{\gamma}_{\nu}(t)
\end{equation}
 is the Student-t density (\ref{eq:student}) with $2\nu$ degrees of
 freedom. The corresponding probability measure is denoted $\sigma_\nu$.
 This formula says that $\sigma_\nu$ is
subordinated to the Gaussian semigroup by an inverse Gamma distribution,
and it implies the infinite divisibility of Student-t from the infinite
divisibility of inverse Gamma. Since the Laplace transformation is
one-to-one, it is clear that if two probabilities $\gamma_1,\gamma_2$
on $]0,\infty[$ lead to the same subordinated density
$$
\int_0^\infty g_t(x)\,d\gamma_1(t)=\int_0^\infty g_t(x)\,d\gamma_2(t),
\quad x\in\mathbb R,$$ 
then $\gamma_1=\gamma_2$.

If we denote $\tau_{a}(x)=ax$, the distribution
$\tau_{a}(\sigma_{n+\tfrac12})*
\tau_{1-a}(\sigma_{m+\tfrac12})$
is given in (\ref{eq:conv}). 
However note that $\tau_{a}(g_{t}(x)dx)=g_{ta^{2}}(x)dx$ so 
\begin{equation}\label{eq:rescale}
\tau_{a}(\sigma_{\nu})=\int_{0}^{\infty}g_{ta^{2}}(x)d\tilde{\gamma}_{\nu}(t)\;
dx,
\end{equation}
 hence \[
\tau_{a}(\sigma_{\nu_1})*\tau_{1-a}(\sigma_{\nu_2})=
\int_{0}^{\infty}\int_{0}^{\infty}(g_{ta^{2}}dx)*(g_{s(1-a)^{2}}dx)\,
d\tilde{\gamma}_{\nu_1}(t)d\tilde{\gamma}_{\nu_2}(s)\]
 \[
=\int_{0}^{\infty}\int_{0}^{\infty}(g_{ta^{2}+s(1-a)^{2}}dx)\,d\tilde{
\gamma}_{\nu_1}(t)d\tilde{\gamma}_{\nu_2}(s)
\]
\[=\int_{0}^{\infty}g_{u}(x)\;
 d\tau_{a^{2}}(\tilde{\gamma}_{\nu_1})*\tau_{(1-a)^{2}}(\tilde{\gamma}_{\nu_2})(u)\; dx.\]
 Therefore, using (\ref{eq:rescale}) we see that for
 $\nu_1=n+\tfrac12,\nu_2=m+\tfrac12$ with $n,m=0,1,\ldots$ the formula
 (\ref{eq:conv}) rewritten as
 \[
\tau_{a}(\sigma_{n+\tfrac12})*\tau_{1-a}(\sigma_{m+\tfrac12})=\sum_{k=n\wedge m}^{n+m}
\beta_{k}^{(n,m)}(a)\sigma_{k+\tfrac12}\]
 is equivalent to 
\begin{equation}\label{eq:conv1}
\tau_{a^{2}}(\tilde{\gamma}_{n+\tfrac12})*\tau_{(1-a)^{2}}
(\tilde{\gamma}_{m+\tfrac12})=\sum_{k=n\wedge m}^{n+m}
\beta_{k}^{(n,m)}(a)\tilde{\gamma}_{k+\tfrac12}.
\end{equation}
 This shows that Theorem \ref{thm:3} is equivalent to the following
 result about inverse Gamma distributions:

{\it The distribution of $a^2Z_n+(1-a)^2Z_m$, where $Z_n,Z_m$ are
  independent inverse Gamma random variables with distribution 
(\ref{eq:inversegaussian}) for $\nu=n+\tfrac12,m+\tfrac12$
respectively, has a density which is a convex combination
of inverse Gamma densities.}

\medskip
This result can be extended to the multivariate Student-t distributions
as follows. A rotation invariant $N-$variate Student-t probability
density is given
for $\bold{x}=(x_1,\ldots,x_N)\in\mathbb R^N$ by
\[
f_{N,\nu}\left(\bold{x}\right)=A_{N,\nu}
\left(1+|\bold{x}|^{2}\right)^{-\nu-\frac{N}{2}},\quad 
A_{N,\nu}=\frac{\Gamma\left(\nu+\frac{N}{2}\right)}
{\Gamma\left(\nu\right)(\Gamma(\tfrac{1}{2}))^N},\]
where 
\[
\langle\bold{x},\bold{y}\rangle=\sum_{i=1}^N x_iy_i,\quad |\bold{x}|=
\left(\langle\bold{x},\bold{x}\rangle\right)^{\tfrac12}, 
\quad\bold{x},\bold{y}\in \mathbb R^N.
\]
It is easy to verify that $f_{N,\nu}(\bold{x})$ is subordinated to the
N-variate Gaussian semigroup
$$
g_{N,t}(\bold{x})=(4\pi t)^{-\tfrac{N}{2}}\exp\left(-\frac
{|\bold{x}|^2}{4t}\right),\quad t>0,\bold{x}\in\mathbb R^N
$$
by the inverse Gamma density (\ref{eq:inversegaussian}), i.e.
\begin{equation}\label{eq:subordination1}
f_{N,\nu}(\bold{x})=\int_{0}^{\infty}g_{N,t}(\bold{x})\,
 d\tilde{\gamma}_{\nu}(t).
\end{equation}
Therefore the characteristic function is given by
\begin{equation}\label{eq:multi}
\int_{\mathbb{R}^{N}}e^{i\langle\bold{x},\bold{y}\rangle} f_{N,\nu}\left(\bold{x}\right)d\bold{x}=k_\nu(|\bold{y}|)
\end{equation}
generalizing (\ref{eq:char}). In fact
$$
\int_{\mathbb{R}^{N}}e^{i\langle\bold{x},\bold{y}\rangle}
f_{N,\nu}\left(\bold{x}\right)d\bold{x}=\int_0^\infty\left(\int_{\mathbb{R}^{N}}e^{i\langle\bold{x},\bold{y}\rangle}
  g_{N,t}\left(\bold{x}\right)\,d\bold{x}\right)\,
d\tilde{\gamma}_\nu(t)
$$
$$
=\int_0^\infty e^{-t|\bold{y}|^2}\,d\tilde{\gamma}_\nu(t)
$$
and the result follows by (\ref{eq:GR}).

As a conclusion, the Theorems 2.1, 2.2 and 2.4 apply in the multivariate
case. For example, an equivalent form of (\ref{eq:conv}) writes as follows:
with $0<a<1$,\[
\frac{1}{a^{N}}f_{N,n+\frac{1}{2}}\left(a^{-1}\bold{x}\right)*\frac{1}{\left(1-a\right)^{N}}f_{N,m+\frac{1}{2}}\left(\left(1-a\right)^{-1}\bold{x}\right)=
\sum_{k=n\wedge
  m}^{n+m}\beta_{k}^{\left(n,m\right)}\left(a\right)f_{N,k+\frac{1}{2}}
\left(\bold{x}\right).\]

\author{Department of Mathematics, University of Copenhagen,
  Universitetsparken 5, DK-2100, Copenhagen, Denmark.}\\
\noindent
\author{Email: berg@math.ku.dk}
 \bigskip

\noindent
\author{Laboratoire d'Informatique IGM, UMR 8049, Universit\'e de
 Marne-la-Vall\'ee, 5 Bd. Descartes, F-77454 Marne-la-Vall\'ee Cedex,
 France.}\\
\noindent
\author{Email: vignat@univ-mlv.fr}
\end{document}